\input colordvi       

\documentclass[amstex,a4]{article}
\usepackage{amsfonts}
\usepackage{epsfig}
\usepackage{floatflt}
\usepackage{color}   
\usepackage{graphicx}
\newtheorem{lemma} {\framebox{\it Lemma}}
\newtheorem{theorem} {\framebox{\it Theorem}}
\newcounter{equacaoUm}

\newcounter{equacaoTres}
\newcounter{equacaoQuatro}
\newcounter{equacaoCinco} 
\newcounter{equacaoSeis}  

\newcounter{equacaoSete}
\newcounter{equacaoOito}

\begin{document}

\pagestyle{empty}
\nocite{Auban:Pierre}
\nocite{DahmenLatour}
\nocite{Calc}
\nocite{Praciano:Um}
\nocite{Rudin}
\nocite{Schumaker}
\nocite{GPL}
\nocite{Gnuplot}

\title{Convolutions power of a characteristic function}
\author{Neves, A. J.\thanks{U Aveiro - Portugal jorgeneves@ua.pt} \\
Praciano-Pereira, T.\thanks{U Vale do Acara\'u - Sobral - Brazil - tarcisio@member.ams.org}}
\date{\today}
\maketitle

\begin{abstract}

This paper deals with the 
convolution powers of the characteristic function of $[0,1], \chi_{[0,1]}$ and its 
function-derivatives. The importance
that such convolution products have can be seen, for an instance, at
\cite{DahmenLatour} where there is the need to find the best differentiable
splines as regularization tool to be used to expand functions or in the
spectral analysis of signals. Another simple application of these
convolution powers is in the construction of a
partition of the unity of a very  high class of differentiability. Here, 
some properties of these convolution powers have been obtained  which easily led
 to write the algorithm to produce convolution powers of $\chi_{[0,1]}$ and their 
function-derivatives. This algorithm was  written in {\tt calc}, the program
 published under GPL that is free to  download.

\par\underline{\bf key words:} \ {\it 
convolution power, b-splines,compact support, GPL
}.

\end{abstract}


\pagestyle{headings}

\section[Notation and conventions]{Introduction} 
\label{Intro}
\par \hskip 0.5cm The $n$-th convolution power of the characteristic function of
an interval is an $(n-1)$-b-splines, this is the subject of this paper.
This first section will be a general description of the work done, the notation used, finishing with three lemmas and a theorem, which
are the main tools used in this work. 

Then we proceed in section \ref{Powers} with the calculus of the 
powers 2,3,4 in a way to uncover the recursion which will be finally clear in the
calculus of the third power. The forth power will be a
preparation for the final theorem (algorithm) which will be presented at
the section \ref{NPower}. The section \ref{Powers} is, thus, 
preparatory for the proof 
at section \ref{NPower}, and  was thought
to permit the reader to construct the algorithm for himself, 
so the reader can skip to the  section
\ref{NPower}, if he find that this subject is familiar and 
if he want only to see the main result.

In the section \ref{Program} it will be presented a
computer program which is a machine translation of the theorem
stated at section \ref{NPower}. Finally, the last section, contains some  applications
 suggested for the material in this paper and a characterization of a general space of univariated splines.

 We have used {\tt calc} 
\cite{Calc}, a computer language specialized to construct mathematical algorithms with
{\em infinite precision} this meaning that it can deals with arbitrarily
long integers which is fundamental to this work as it will be necessary
to use factorial for any value of $n$ in the calculus of the convolution powers.
In section \ref{Program} you can
see a graphic output of the program for $n=30$.

\subsection[The history]{Some history}
\label{History}
\par \hskip 0.5cm Splines are pieces of polynomials pasted together in
conditions to make a continuous, and even more, a highly differentiable
function.
In the last section you may find a definition for a class
of splines.
 It is a simple idea which turned into a very powerful tool in such
diverse areas as {\em geometric design} or {\em approximation theory}. 
In the background
there is a famous theorem of Weierstrass stating that any continuous
function can be arbitrarily approximated by polynomials. There was a huge work to
make this simple and easy, and we can consider a bad example the Taylor
polynomials, which make a very good approximation of a function in a sole
point but which is good enough to be used to approximate periodic functions
like sine. The bad thing with {\em classical polynomial approximation} is that computers
do not like it that much because it needs a very high degree to get
something useful. You can remember Lagrange polynomial
interpolation which does a good job if the points to be linked together are
not jumping to much and the critical situation is at the boundary of the
interval which defines the region where the points have been selected. 
This is a particular problem to {\em classical polynomial approximation} 
because just around
the boundary of the region containing its roots a polynomial will {\em behave
 wildly}.  If
you have selected $n$ points the corresponding Lagrange polynomial is of degree $n-1$.

Then splines came out to be the solution. With some pieces of first degree 
polynomials you may have a very good approximation of data collected 
from a sensor, provided that you are only interested in the total amount
of the phenomenon, its integral. Well this was the trick, paste together
several piece of polynomials instead of a unique monolithic polynomial
as a Lagrange polynomials. The piecewise first degree polynomials 
are 1-splines, known as polygonal
lines.

Switching to polynomials of second degree 
we can already have a better approximation and the famous algorithm of Runge-Kutta
does exactly this for solving differential equations, though not 
with splines. 
But degree two
polynomials have a big lack of flexibility since they do not provide a change
of convexity in each piece, they are 2-splines.  
So, the next step is to paste together
pieces of polynomials of third degree and this has proved to be a good
solution for a large number of problems needing polynomial approximation: 3-splines.  

At this
point it should pop up a definition for this little thing called splines
that started to be used in the seventies and today is a very big network of
research in Mathematics. To be an $n$-splines, the pieces of polynomial 
should be of degree less or equal to $n$,
there can be a straight line infiltrated in the
business. The second condition is that the resulting function has
to be of class $C^{n-1}$.

It is simple like that but the main
literature about splines, for example, the historic book of Larry Schumaker
\cite{Schumaker}, is not easy to be read though it is at the beginning
of the history. This is the point with this paper, it is starting by the end,
not from the beginning, and there are some powerful tools that had been added
to the history to make splines easy.

\subsection[Contribution]{The contribution of this paper}
\label{Contribution}

\par \hskip 0.5cm In chapter 6 of \cite{Auban:Pierre}, the author has shown that 
the n-th convolution power of the characteristic function of the interval
$[0,1]$ is an $(n-1)$-splines having as support the interval $[0,n]$ and he
starts by showing that
 the convolution of two characteristic functions results in a
1-splines. Here you are the basic tool to build splines: the 
convolution. Convolution is a ``multiplication'', it satisfies the
distributive law relatively to the sum, hence a sum of characteristic
functions of intervals can be multiplied by a \textit{selected characteristic
function of an interval}, using convolution, to produce a 1-splines. Iterating this process you will reach
an n-splines for the desired $n$. But this take a lot of time and computer
power and the alternative is to construct a {\em splines-bases} to generate
an appropriate {\em space of splines} and this is the tool that this paper is going to 
produce: an approximate unity to build up bases for vector spaces of splines.
This is an well known application  that will be only mentioned at the
last section as complement to the main result of this paper.
These convolution powers can be easily be transformed into
a {\em sharp approximate unity} with respect to the convolution product and this
will be mentioned again at  section \ref{Conclusions}.

Here we have a slightly different methodology from the one 
at \cite{Auban:Pierre},  besides presenting
four properties of $f=(\chi_{[0,1]})^{n}$, which are the lemmas
and the theorem at section \ref{Powers},
 the work is directed to  a 
recursive program that is able to produce the whole formula not 
only of $f=(\chi_{[0,1]})^{n}$ but for all its function-derivatives,
having as input the power $n$ which is needed. 
This is the main result of the paper, theorem \ref{TeoremaQuatro}, in the  section \ref{NPower} followed
by the correspondent algorithm, a computer program,  in section  \ref{Program}. The program
produces 
functions for {\tt gnuplot} \cite{Gnuplot}, 
which may be run to show the graphics of $f$  and is able to write the matrix of  the coefficients of $f$ and all its function-derivatives.

There are recent works, see \cite{DahmenLatour} for example, which show the importance of convolutions in ``spectral analysis''
mainly in the context of wavelets. 
This is an example of the 
importance of having a simple way to present convolution powers and 
a good algorithm to produce them. This paper offers that contribution
in this matter in the particular case of convolution of characteristic
functions of intervals of the line. 

The whole construction is 
done step by step, in a way to permit the reader to construct by himself
the method. 

Our main motivation with this paper is the construction of 
a tool to be used in the {\em analysis of differential operators} and in the {\em spectral analysis of wavelet}
and we think this is underway.

\section[Convolutions power of a characteristic function] 
		{The convolution powers of $\chi$ } 
\label{Powers}

\par \hskip 0.5cm  In this section, first, we introduce the notation and basic results necessary to 
proof the main theorem of the paper, following this we are going to calculate the 
powers 2, 3, 4, in order to formulate the general result.

\subsection[Notation]{Basic results}
\label{Notation}

\par \hskip 0.5cm Here you are the basic notation that will  be stated 
together with three Lemmas which
are the underground tools for the theorems of the paper.

The characteristic function of the interval $[0,1]$ 
will be refereed to using one
of the following notations: 
\begin{equation}
\chi_{0} = \chi_{[0,1]} = \chi
\end{equation} 
 
It will be used the index
when necessary to indicate a translation of $a \in \mathbf{R}$, then
\begin{equation}
\chi_{a}(x) = \chi_{0}(x-a) = \chi(x-a); 
\end{equation} 

The use power notation for ``convolution power'' has no risk of 
 possible confusion with the ``arithmetic power'', since the
context will indicate what is the meaning in the few occasions in which the ``arithmetic power'' will be used. Hence
$\chi = \chi^{1}$ is the first power, as usual, and we can add a
sense to  $\chi^{0}$ as the constant function. 

The symbol
$f = \chi^{n}$ will be used to simplify the notation of the $n$-th convolution
power, which is the main goal of this paper, hence $f$ do represent the $n$-th convolution power of $\chi$, 
consistently, in the whole paper.

The second convolution power of the characteristic function is the
{\em triangle
function}, a 1-splines, 
whose derivative is the the difference of two the translations
of the characteristic
function: $\chi_{0}-\chi_{1}$. This is a consequence of the following
lemmas.

\begin{lemma}[Convolution and derivative]
\label{LemmaUm}

$$ \frac{d}{dx} (f*g) = \frac{d f}{dx}*g = f*\frac{d g}{dx}
$$
\end{lemma}

\begin{lemma}[Distribution derivative of characteristic function]
\label{LemmaDois}

$$\frac{d \chi}{dx} = \delta_{0} - \delta_{1}
$$
\end{lemma}

\begin{lemma}[Convolution with $\delta_{a}$ and translation]
\label{LemmaTres}

$$
	(\delta_{a}*f)(x) = f(x-a)
$$
\end{lemma}

Using these lemmas it can be proved:

\begin{theorem}[Derivative of the second convolution power]

\label{TeoremaUm}

If 
$$f = \chi*\chi = \chi^{2}$$ 
then 
$$f' = \chi_{0}-\chi_{1}$$

\end{theorem}

In lemma \ref{LemmaTres}, when the translation parameter is $a=0$, 
we can see that the {\em Dirac delta} measure is the unity of convolution
as an element of an extension ring of the convolution ring  of functions
which does not have  a unity.

\subsection[The second power]{The second convolution power}
\label{SecondPower}

\par \hskip 0.5cm These three properties stated as lemmas 
produce an easy way to get the second 
convolution power of  $\chi$, the above mentioned triangle function,
(compare with the heavy work at \cite[chapter 6]{Auban:Pierre}. Moreover, this method will be 
used all along the paper, repeatedly, the main reason being that we shall
have a succession of integrals whose domains of integration are disjoints, they
coincide, in fact, exactly on a single point.   
\begin{eqnarray}   
\setcounter{equacaoUm}{\arabic{equation}}  
& f' = (\chi*\chi)^{'} = \chi*\chi^{'}= \chi*(\delta_{0} - \delta_{1} )  \\
& f' = (\chi*\chi)^{'} = \chi*\delta_{0} -  \chi_{0}*\delta_{1}  \\
& f' = (\chi*\chi)^{'} =  \chi -  \chi_{1} \\
& f(x) = \chi*\chi(x) = \chi^{2}(x) = 
	\int\limits_{0}^{x} \left( \chi(t) - \chi_{1}(t)\right) dt  \\
&f(x) =\chi^{2}(x)  = 
	\left\{ \begin{array}{ll}
		x<0 & 0 \\
		x \in [0,1] &  \int\limits_{0}^{x} dt  \\
		x \in [1,2]  & a - \int\limits_{1}^{x}\chi_{1}(t) dt \\
		x > 2 & 0 \\
		\end{array} \right.\\
& f(x) =\chi^{2}(x)  =
	\left\{ \begin{array}{ll}
		x<0 & 0 \\
		x \in [0,1] &  x  \\
		x \in [1,2]  & 1 - \int\limits_{0}^{x-1} dt \\
		x > 2 & 0 \\
		\end{array} \right.\\
&f(x) =\chi^{2}(x) =
	\left\{ \begin{array}{ll}
		x<0 & 0 \\
		x \in [0,1] &  x  \\
		x \in [1,2]  & 1 - (x-1)  \\
		x > 2 & 0 \\
		\end{array} \right.
\end{eqnarray} 
and now it is easy to see that $\chi^{2}$ is the so called triangle function
centered at 1 with support $[0,2]$ and height $1$ at 1. You may observe
the presence of the symbol ``$a$'',
\addtocounter{equacaoUm}{4} at equation (\arabic{equacaoUm}),
\addtocounter{equacaoUm}{-4}
 which is the value of the 
integral on $[0,1]$ and is the initial condition of the integral
on $[1,2]$. This trend
will be repeated further and further for all convolution powers.
The symbols ``$a,b,c, \ldots$'', will be used temporally, before  
their values can be expressed. 

The triangle function can be defined as 
$
\chi^{2}(x) =x\chi + (2 - x)\chi_{1}
$.
This expression is not useful in computations, but it will be of good
help in the mathematical formulation, in the logic preparation which will 
produce the algorithm. 

Observe that the (distribution) derivative of  $\chi$ is $\delta_{0} - \delta_{1}$
where $\delta_{a}$ is the Dirac measure concentrated at the point $a$.
The fact that the ``convolution with the Dirac measure is a translation of
the other convolution factor'' has been used in equation
\addtocounter{equacaoUm}{2}  (\arabic{equacaoUm}).
\addtocounter{equacaoUm}{-2} 

An easy and beautiful theorem, extending theorem \ref{TeoremaUm}, can be  stated, immediately.

\begin{theorem}[Derivative of $f$]{the $n$-th convolution power}
\label{TeoremaDois}
\begin{equation}
\frac{d}{dx}f(x) = \frac{d}{dx}\chi^{n}(x) = \chi^{n-1}(x) - \chi_{1}^{n-1}(x)
\end{equation}
\end{theorem}
the {\em derivative of a convolution power $f$ of $\chi$ is the difference of 
translates of the previous power}. The proof is a direct
application of  Lemmas \ref{LemmaUm}, \ref{LemmaDois} and \ref{LemmaTres}. 
This theorem  will be used several times in the sequence.

A notation will be introduced to simplify the expression
of $f$. Put $P_{1}(x)=x$ 
then the equation
\addtocounter{equacaoUm}{6}(\arabic{equacaoUm})
\addtocounter{equacaoUm}{-6} reads:
\begin{equation}
f(x) = \chi^{2}(x) =
          \left\{ \begin{array}{ll}
		x<0 & 0 \\
		x \in [0,1] &  P_{1}(x)  \\
		x \in [1,2]  & 1 - P_{1}(x-1)  \\
		x > 2 & 0 \\
                    \end{array} \right.
\end{equation}

Observe that the second convolution power  has been defined  in terms of 
a linear combination of a polynomial of degree one and its translates (the plural will be valid in a near future\dots), 
this highlights that $f$,
the second convolution power, is 1-b-splines as particular
case of the first sentence of the paper, that ``$n$-th convolution
power is an $(n-1)$-b-splines''. It will be possible to go further
using this kind of polynomials and the same will be done  
with the next power
so to make clear the recursion of the process.

\subsection[The third power]{The third convolution power}
\label{ThirdPower}

\par \hskip 0.5cm Before going ahead let us stress that the whole 
secret is contained
in first calculating the derivatives of the next power which we already have,
thinking that, to calculate the next power first we must have the 
previous one\dots 
So let's 
calculate the derivatives of the third convolution power to have:

\begin{eqnarray}
\setcounter{equacaoTres}{\arabic{equation}}  
& \frac{d f}{dx} = \frac{d}{dx}\chi^{3}(x) =  \chi^{2}(x) -
\chi^{2}(x-1) \\
& \frac{d^{2} f}{dx^{2}} = \frac{d^{2}}{dx^{2}}\chi^{3}(x) = 
\chi(x) - 2\chi(x-1) + \chi(x-2);
\end{eqnarray}
where we can see the binomial coefficients: 1,-2,1 as usual when
calculating derivatives of products and, here, there is a ``difference''
in the ``game''. 

This suggests, by the way,  another beautiful and easy
theorem:

\begin{theorem}[The $n-1$ derivative]{of the $n$ convolution power}
\label{TeoremaTres}
\begin{enumerate}
\item  The $n$-th convolution power, $f=\chi^{n}$, is an $(n-1)$-b-splines of
compact support $[0,n]$, hence $f$ has $n-1$ function-derivatives of whose
$n-2$ are continuous.

\item The derivative of order $n-1$ of the $n$-th convolution power $f$, 
which
is not continuous, is given by the linear combination of translates
of $\chi$
\[
\frac{d^{(n-1)} f}{dx^{(n-1)}}   = \sum\limits_{k=0}^{n-1} (-1)^{k}(^{n-1}_{k})\chi_{k}
\]
\end{enumerate}
The coefficients of this linear combination are the elements of the line
of order $n-1$ of the ``alternate Pascal Triangle'', these elements
are the coefficients of $(x-1)^{n-1}$, (the line of order zero
has only one element, 1).
\end{theorem}

The proof comes by recursion with repeated application of 
the lemmas \ref{LemmaUm},  \ref{LemmaDois}, \ref{LemmaTres}
and theorem \ref{TeoremaUm}.

The calculus of the iterated integral of $f''$  will be easy
as it is again a succession of integrals
whose integration domains are disjoints (coincide on a set of measure zero),
one point set.

Hence the third convolution power is twice the integral of its 
second derivative
starting at the initial condition $x=0$. We know in advance
that the support of this convolution is $[0,3]$, from \cite{Auban:Pierre}, 
and that there
will be different equations each time we go over an integer boundary.
The initial condition of each integral is the {\em total value} of the previous integral to make a continuous function.
This will be represented by the symbols $a,b$, temporally.

In condition to obtain the recursion formula, it is necessary to
introduce here a {\em rather complicated notation}, but intuitive, indeed. 
It will protect the
parts of the expression which will pop up. Observe that the second
derivative of the third convolution power is not continuous as 
a linear combination of translations of $\chi$ with coefficients
being  the {\em alternate binomial
numbers}, as in the powers of (x-1).  So let the following numbers be
defined:

\begin{equation}
\left\{ \begin{array}{l}
Bin(n,k) = (-1)^{k}(^{n}_{k}) \\
	a_{00}=Bin(2,0); \\ a_{01}=Bin(2,1); \\ a_{02}=Bin(2,2);\\
	\end{array} \right.
\end{equation} 

With this notation it may be written

\begin{equation}
\frac{d^{2} f}{dx^{2}} = 
\frac{d^{2} }{dx^{2}}\chi^{3} =  f''(x) = \left\{ \begin{array}{ll}
	x <0      & 0 \\
	x \leq 1  &  a_{00}\chi(x)\\
	x \leq 2  &  a_{01}\chi(x-1)\\
	x \leq 3  &  a_{02}\chi(x-2)\\
	x > 3  &  0 \\ \end{array} \right.
\end{equation}

Integrating twice this second derivative
we shall obtain
 the
formulation of $f$. At this point, 
$f$ is the second convolution
power of $\chi$. A new set of 
coefficients will be defined in the run and, 
for a while, the convention, 
state previously, to use the symbols $a,b$ as the initial 
conditions for the second and third integrals, will be used.

\begin{eqnarray}
\setcounter{equacaoCinco}{\arabic{equation}}  
& f'(x) = 
\int_{0}^{x}\left( a_{00} \chi(t) + a_{01} \chi(t-1) + a_{02}\chi(t-2) \right) dt \\
& f'(x) = \left\{ \begin{array}{ll}
	x <0  	& 0 \\
	x \leq 1	&  a_{00}x \\
	x \leq 2  &   a + a_{01}\int_{1}^{x}  \chi(t-1) dt \\
	x \leq 3  &   b + a_{02}\int_{2}^{x} \chi(t-2) dt \\
	x > 3    &  0 \\
	\end{array}\right. \\ 
& f'(x) = \left\{ \begin{array}{ll}
	x <0  	& 0 \\
	x \leq 1	&  a_{00}x \\
	x \leq 2  &   a + a_{01}\int_{0}^{x-1}  \chi(t) dt \\
	x \leq 3  &   b + a_{02}\int_{0}^{x-2} \chi(t) dt \\
	x > 3    &  0 \\
	\end{array}\right. \\ 
& f'(x) = \left\{ \begin{array}{ll}
	x <0  	& 0 \\
	x \leq 1	&  a_{00}x \\
	x \leq 2  &   a_{00}  +  a_{01}(x-1)  \\
	x \leq 3  &   a_{00} + a_{01}   + a_{02}(x-2)  \\
	x > 3    &  0 \\
	\end{array}\right. \\ 
& f'(x) = \left\{ \begin{array}{ll}
	x <0  	& 0 \\
	x \leq 1  &  a_{10}x \\
	x \leq 2  &   a_{11}  + a_{01}(x-1) \\
	x \leq 3  &   a_{12}  + a_{02}(x-2)  \\
	x > 3    &  0 \\
	\end{array}\right. \\ 
& \left\{ \begin{array}{ll} 
	a_{10} = & a_{00}; \\
	a_{11} = & a_{10}P_{1}(1); \\
 	a_{12} = & a_{11} + a_{01}P_{1}(1);\\ \end{array} \right.
\end{eqnarray}

The equation\addtocounter{equacaoCinco}{5} (\arabic{equacaoCinco})
\addtocounter{equacaoCinco}{-5}
 can be written a bit different in a way to show that
the sum of certain indexes is constant in each row.
This rule is valid with {\em two exceptions}
\footnote{though that, at this point, the number of exceptions is very
high relatively to the number of cases, 1, for which the rule is in force \dots}
 which are meaningless to the
definition of the algorithm:
\begin{enumerate}
\item The first row of the matrix $a_{0j} = Bin(n,j)$;
\item The two first elements in in each row
$$
	a_{i0} = a_{00}; a_{i1}=a_{i0}=a_{00}
$$
\end{enumerate}

\begin{equation}
\setcounter{equacaoSeis}{\arabic{equation}}  
\left\{ \begin{array}{ll} 
	a_{10} = & a_{00}P_{1}(1);\\ 
	a_{11} = & a_{10}P_{0}(1);\\ 
	a_{12} = & a_{11}P_{0}(1) + a_{01}P_{1}(1); \\  \end{array} \right.
\end{equation}
In equation (\arabic{equacaoSeis}) it is being used the family of polynomials
\begin{equation}
\setcounter{equacaoSeis}{\arabic{equation}}  
P_{k}(x) = x^{k}/k!;
\end{equation}

In the equations to define the coefficients $a_{ij}$ 
we have in each row a sum of products $a_{ij}P_{k}(1)$ and the sum 
$i+k$ is the order of the row when $j\geq 2$ (the above mentioned exception).

To obtain a
new row 
 $a_{ij}$, we have to  use all the possibilities in ``$i+k =$order of the row''.
This law will be valid in the sequence.

Calculating the next integral we shall have:

\begin{eqnarray}
\setcounter{equacaoSete}{\arabic{equation}}  
& f(x) = \left\{ \begin{array}{ll}
	x < 0  	& 0 \\
	x \leq 1	&  a_{10}\int\limits_{0}^{x} t dt \\
	x \leq 2 & a + a_{11}\int\limits_{1}^{x}dt + a_{01}\int\limits_{1}^{x}(t-1) dt \\
	x \leq 3 & b + a_{12}\int\limits_{2}^{x}dt + a_{02}\int\limits_{2}^{x}(t-2) dt\\
	x > 3    &  0 \\
	\end{array}\right. \\
& f(x) = \left\{ \begin{array}{ll}
	x < 0  	& 0 \\
	x \leq 1	&  a_{10}P_{2}(x) \\
	x \leq 2 & a_{10}P_{2}(1) + a_{11}\int\limits_{0}^{x-1}dt 
			+ a_{01}\int\limits_{0}^{x-1}t dt \\
	x \leq 3 & b + a_{12}\int\limits_{0}^{x-2}dt + a_{02}\int\limits_{0}^{x-2}t dt\\
	x > 3    &  0 \\
	\end{array}\right. \\
& f(x) = \left\{ \begin{array}{ll}
	x < 0  	& 0 \\
	x \leq 1	&  a_{20}P_{2}(x) \\
	x \leq 2 & a_{20}P_{2}(1) + a_{11}P_{1}(x-1) + a_{01}P_{2}(x-1) \\
	x \leq 3 & b + a_{12}P_{1}(x-2) + a_{02}P_{2}(x-2)\\
	x > 3    &  0 \\
	\end{array}\right. \\
& f(x) = \left\{ \begin{array}{ll}
	x < 0  	& 0 \\
	x \leq 1	&  a_{20}P_{2}(x) \\
	x \leq 2 & a_{20}P_{2}(1) + a_{11}P_{1}(x-1) + a_{01}P_{2}(x-1) \\
	x \leq 3 & a_{20}P_{2}(1) + a_{11}P_{1}(1)   + a_{01}P_{2}(1) 
			+ a_{12}P_{1}(x-2) + a_{02}P_{2}(x-2)\\
	x > 3    &  0 \\
	\end{array}\right. \\
& f(x) = \left\{ \begin{array}{ll}
	x < 0  	& 0 \\
	x \leq 1	&  a_{20}P_{2}(x) \\
	x \leq 2 & a_{21}P_{2}(1) + a_{11}P_{1}(x-1) + a_{01}P_{2}(x-1) \\
	x \leq 3 & a_{22} + a_{12}P_{1}(x-2) + a_{02}P_{2}(x-2)\\
	x > 3    &  0 \\
	\end{array}\right. \\
\end{eqnarray}

Some remarks here will make things clear in the next step. In equation (\arabic{equacaoSete})
the symbols $a,b$ have been used as explained before, and in equation 
\addtocounter{equacaoSete}{1} (\arabic{equacaoSete}) \addtocounter{equacaoSete}{-1}
 the limits of integration have been changed accordingly with the interval of definition,
and the liberty of translation invariance of Lebesgue measure has been used
to change the limits of  $\int\limits_{0}^{x-2}dt$ in conditions to have the same
setting in both integrals. This leads to a translation of the
polynomial $P_{i}$, naturally. This will be further commented in a moment. 
The symbol $a$ has been evaluated as $f(1)$
and the symbol $b$ has been evaluated as
$$ 
b = f(2) = a_{21}P_{2}(1) + a_{11}P_{1}(1) + a_{01}P_{2}(1) 
$$ 
the value 
of the primitive over $[1,2]$ evaluated at point $x=2$. 
On the interval $[2,3]$ it can be written $a_{22}P_{0}(x-2)$
to enhance the property 
$$
a_{mj}P_{i}(x-k); m+i=j
$$ 
where $j$ is the
order of the row.
The new set of coefficients will be defined explicitly:  
\begin{equation}
\setcounter{equacaoOito}{\arabic{equation}}  
\left\{ \begin{array}{ll}
a_{20} & = a_{10}= a_{00}; \\ 
a_{21} & = a_{20}P_{2}(1); \\ 
a_{22} & = a_{21}P_{0}(1) + a_{11}P_{1}(1) + a_{01}P_{2}(1);\\
	\end{array} \right.
\end{equation}
and remember the exception to the rule that the first two coefficients in each row 
have a particular law. With these definitions $f$ is

\begin{equation}
f(x) = \left\{ \begin{array}{ll}
	x <0      & 0 \\
	x \leq 1  & a_{20}P_{2}(x); \\
	x \leq 2  & a_{21}P_{0}(x-1) + a_{11}P_{1}(x-1) + a_{01}P_{2}(x-1)\\
	x \leq 3  & a_{22}P_{0}(x-2) + a_{12}P_{1}(x-2) + a_{02}P_{2}(x-2)\\
	x > 3    &  0 \\ \end{array}\right.\\
\end{equation}
a linear combination of the elements of last row of the matrix $a[i,j]$
with translations of all the polynomials 
$P_{k}; k=0,\ldots,n$, where
$n+1$ is the desired convolution power. The previous raws of this matrix give
the derivatives of $f$ but the summations are being restricted of one
unity in each further derivative. This will be clear in the last section, but
you can conclude from the above calculations that:

\begin{eqnarray}
& f'(x) = \left\{ \begin{array}{ll}
	x <0      & 0 \\
	x \leq 1  & a_{10}P_{1}(x); \\
	x \leq 2  & a_{11}P_{0}(x-1) + a_{01}P_{1}(x-1) \\
	x \leq 3  & a_{12}P_{0}(x-2) + a_{02}P_{1}(x-2)\\
	x > 3    &  0 \\ \end{array}\right. \\
& f''(x) = \left\{ \begin{array}{ll}
	x <0      & 0 \\
	x \leq 1  & a_{00}P_{0}(x); \\
	x \leq 2  & a_{01}P_{0}(x-1)  \\
	x \leq 3  & a_{02}P_{0}(x-2) \\
	x > 3    &  0 \\ \end{array}\right.
\end{eqnarray}
Remark that $P_{0}$ can be replaced by $\chi$, but this is of no help at all.

\subsection[The forth power]{The forth convolution power}
\label{ForthPower}

\par \hskip 0.5cm From this point the recurrence is established:
\begin{enumerate}
\item The convention for indexes follows the syntax of the language 
${\mathtt{C}}^{++}$ for
which the first positive (integer) index is zero, 
this way in a $n \times n$ matrix the indexes run from $0 \ldots n-1$. 
The matrix $(a_{ij});i,j = 0, \cdots,n-1$ where $n$ is the desired power;
\item to use the entries in  the last row of this matrix to define $f(x)$ as combinations of 
translations of $P_{i}(x),  P_{i}(x-k),$ to the interval $[k,k+1]$. 
\item The $n-1$ derivatives of $f$ are obtained, from the raw $n-2$
to the raw zero with increasing order of derivation and decreasing the
number of summations of powers of $P_{i}$ in the same way as $f$.
\item The starting point is the derivative of order $n-1$ which comes
from the line of order $n-1$ of the {\em alternate Pascal triangle}.
\end{enumerate}

The forth power, $f$, of $\chi$ may, now, be written,
directly, without further integrations as an 
example before writing the algorithm formally.

By Theorem \ref{TeoremaDois}, $f$  has
tree derivatives, two continuous derivatives, 
it is 3-splines. So the algorithm will start from the line of order
3 of the ``alternate'' Pascal triangle:

\begin{eqnarray}   
\setcounter{equacaoQuatro}{\arabic{equation}}  
& a_{0 k}= Bin(3,k); k=0 \dots 3;\\
& \left\{ \begin{array}{ll}
	a_{10} & = a_{00}; \\
 	a_{11} & = a_{10}P_{1}(1);\\
  	a_{12} & = a_{11}P_{0}(1) + a_{01}P_{1}(1);\\
  	a_{13} & = a_{12}P_{0}(1) + a_{02}P_{1}(1); \\
	\end{array} \right. \\
&  \left\{ \begin{array}{ll}
	a_{20} & = a_{00}P_{0}(1) \\
	a_{21} & = a_{20}P_{2}(1); \\
	a_{22} & = a_{21}P_{0}(1) + a_{11}P_{1}(1)+ a_{01}P_{2}(1);\\
	a_{23} & = a_{22}P_{0}(1) + a_{12}P_{1}(1)+ a_{02}P_{2}(1);\\
	\end{array} \right. \\
& \left\{ \begin{array}{ll}
	a_{30} 	& = a_{00}P_{0}(1) \\
 	a_{31}	& = a_{30}P_{3}(1); \\  
 	a_{32}	& = a_{31}P_{0}(1) + a_{21}P_{1}(1) +  a_{11}P_{2}(1) + a_{01}P_{3}(1); \\
	a_{33}	& = a_{32}P_{0}(1) + a_{22}P_{1}(1) +  a_{12}P_{2}(1) + a_{02}P_{3}(1); \\
	\end{array} \right. 
\end{eqnarray} 
and 
\begin{equation}
f(x) = \left\{ \begin{array}{ll}
x<0  &  0 \\
x<1  &  a_{30}P_{3}(x); \\
x<2  &  a_{31}P_{0}(x-1) + a_{21}P_{1}(x-1) + a_{11}P_{2}(x-1) + a_{01}P_{3}(x-1)\\
x<3  &  a_{32}P_{0}(x-2) + a_{22}P_{1}(x-2) + a_{12}P_{2}(x-2) + a_{02}P_{3}(x-2)\\
x<4  &  a_{33}P_{0}(x-3) + a_{23}P_{1}(x-3) + a_{13}P_{2}(x-3) + a_{03}P_{3}(x-3)\\
x>4  &  0 \\ \end{array} \right.
\end{equation}

Observe that by Theorem \ref{TeoremaUm} the derivatives of $f$ can be 
traced back by difference of translations of the previous powers up to
the last function-derivative as a linear combination of translations
of $\chi$ so we are able to right the equations of $f$ and all its 
function-derivatives very easily. Better said, a program can do that for us,
and the program can be download from a link at \cite{Praciano:Um}.

\section[The power $n$]{The $n$-th convolution power}
\label{NPower}

\par \hskip 0.5cm In the next theorem the indexes are separated with commas
to avoid the confusion: ``$n-1,1$'' for ``$n-1 1$''. The same has to be done at the program.

\begin{theorem}[The $n$-convolution power]{of the characteristic $\chi_{[0,1]}$}

\label{TeoremaQuatro}
The $n$-th convolution power of $\chi = \chi_{[0,1]}$, $f =\chi^{n}$, is $n-1$-b-splines
such that 
\begin{eqnarray*} 
& f(x) = \left\{ \begin{array}{ll}
	0                    &  x < 0 \\
	a_{n-1,0}P_{n-1}(x)  &  x \in [0,1] \\
		    \sum\limits_{j=0}^{n-1} a_{n-1,j}P_{j}(x-k);
		    k = \lceil x \rceil; 
	    &  ( 1 \leq x \leq n) \\
	0 &  x > n   \\
	\end{array}\right.
\end{eqnarray*}
where the coefficients $a_{n-1,j}$ are given by the equations:
\begin{eqnarray}
\setcounter{equacaoCinco}{\arabic{equation}}  
& \left[\begin{array}{llllll}
a_{0,0}   & a_{0,1}   & a_{0,2}  & \cdots 	& a_{0,n-1}  \\
a_{1,0}   & a_{1,1}   & a_{1,2}  &  \cdots 	& a_{1,n-1}  \\
	  & 	  	& 	 &\cdots 	& \\
a_{n-1,0} & a_{n-1,1} & a_{n-1,2} &  \cdots     & a_{n-1,n-1}\\
	\end{array}\right] \\
& a_{0,k} = Bin(n-1,k); k=0, \dots, n-1; Bin(n,p) = (-1)^{p}(_{p}^{n-1}) \\
& a_{1,0} = a_{0,0}; \\
& a_{1,1} = a_{1,0}P_{1}(1); \\
& j \geq 2; j \leq n-1;
 a_{1,j} = a_{1,j-1}P_{1}(1); \\
& a_{2,0} = a_{0,0}P_{0}(1); \\
& a_{2,1} = a_{2,0}P_{1}(1); \\
& j \geq 2; j \leq n;
a_{2,j} = a_{2,j-1} + a_{1,j-1}P_{1}(1) + a_{0,j-1}P_{2}(1); \\
& a_{3,0} = a_{0,0}P_{0}(1); \\
& a_{3,1} = a_{3,0}P_{1}(1); \\
& j \geq 2; j \leq n;
a_{3,j} =   \sum\limits_{i+k=3}   b{i,j-1}P_{k}(1) \\
&  \cdots \\
& a_{n-1,0} = a_{0,0}P_{0}(1); \\
& a_{n-1,1} = a_{n-1,0}P_{1}(1); \\
& j \geq 2; j \leq n-1;
a_{n-1,j} =   \sum\limits_{i+k=n-1}   a_{i,j-1}P_{k}(1)
\end{eqnarray} 

\underline{\bf The proof:}

By theorem \ref{TeoremaDois}, $n$-th convolution power, which is a $(n-1)$-b-splines, has  
$n-1$ function-derivatives, of which the first $n-2$ are continuous, 
and the starting point is line of order $n-1$ of the {\em alternate Pascal triangle}
defined by the numbers $Bin(n-1,k)$. This is a consequence of two facts:
\begin{enumerate}
\item In the combination of lemmas \ref{LemmaUm} and \ref{LemmaDois} we have 
a relation of the type of {\em Stifel rule}, the defining relation of Pascal
Triangle,
modified to produce the alternate numbers $Bin(n-1,k)=(-1)^{k}(_{k}^{n-1})$.
So after $n-1$ derivatives we are at line of order $n-1$ of 
the {\em alternate Pascal triangle}.
\item Each new derivative is a linear combination of the previous convolution
power by  theorem \ref{TeoremaUm} so the order of convolutions power will
decrease to reach order one of derivation, the non-continuous function-derivative.

\end{enumerate}
This explains the first row of the $(n-1) \times (n-1)$-matrix $b[ji]$.
The other rows are obtained by integrating the combination of 
translations of the $P_{j}$ all along of the integers of the support 
$[0,n]$ of $f$, and at each new integration the total variation of
the previous integral is used as the initial value, with $j$ being the
order of the row, $j = n-1 -k$ where $k$ is the order of the derivative.

\end{theorem}

\section[The computer program]{A computer program}
\label{Program}

\par
\begin{figure}[ht]
\center
\includegraphics[width=12cm,height=10cm]{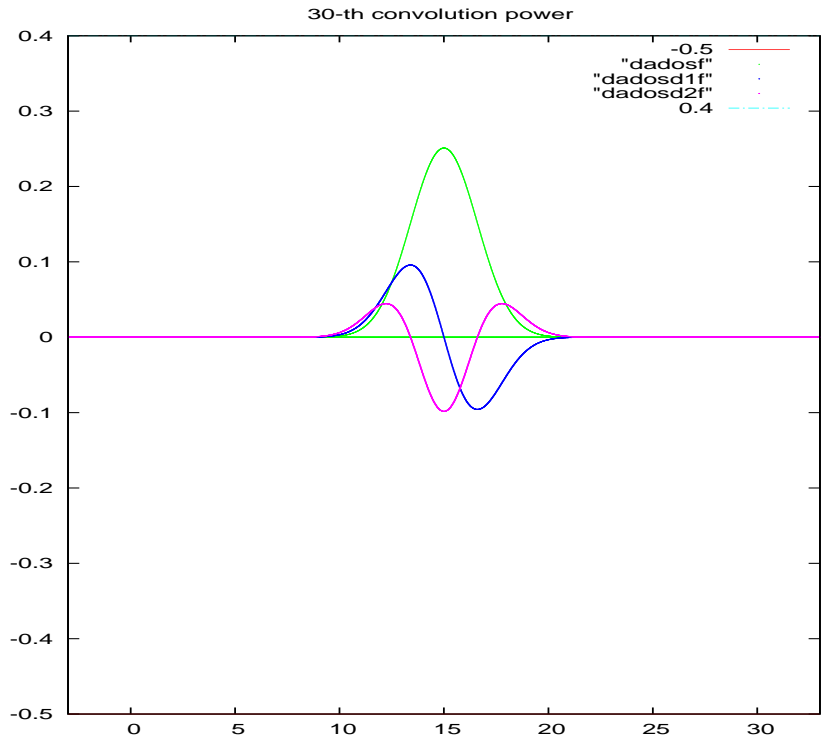}
\caption{\footnotesize graphics of convolution power 30, $f, f' f''$  }
\label{ConvolutionPowersUmUm}
\end{figure}

 \hskip 0.5cm  There is a program written in {\tt Calc}, \cite{Calc},
which can be download from \cite{Praciano:Um}, and here you
are the output of the program at the figure
(\ref{ConvolutionPowersUmUm})    
\leavevmode
for the choice of $n=30$. This is the graphic of $f,f', f''$ when $f=\chi^{30}$.
Time of this job: a few seconds.

To have the graphics of $f,f',f''$, of course, for $n \geq 3$,
do the following:
\begin{enumerate}
\item Start {\tt Calc}, then you shall have a terminal of the {\em calculator}
to make inputs. {\tt Calc} is iterative, it is an {\em interpreted language};
\item Supposing you have the code at the current directory named to 

\par \hskip 2.5cm {\tt convolution\_powers.calc} 

\noindent execute at the terminal of {\tt Calc}:

\par \hskip 2.5cm {\tt read convolution\_powers.calc}

``execute'', means, ``write the above at the terminal and hit enter''.
This will make {\tt Calc} learn all
the functions defined at {\tt convolution\_powers.calc} and you shall be able
to execute any of them. 
There is one function called {\tt main()} which
is written to do a big part of the job so you can have what is promised in
the paper. 

\item Execute ``main(n)'' 
selecting the desired value of $n$  
at the terminal of {\tt Calc}
and you will see $f,f',f''$. This will work {\em out of the box}, under Linux.

You may experience some problems but then read the code and you
shall see a suggestion to go around them.

\item If you want see the matrix, execute:

\begin{itemize}
\item a = cria\_matriz(n); either repeat the the value of n or put a new one.
\item {\tt imprime\_matriz(a,n) } the parameter ``a'' has been created in the
previous step, but you still need to put the same value of $n$ as parameter.
This indicates that the program is not well planned...Perhaps somebody
will produce something better from this draft.
\end{itemize}

\item It is pointless to show all the $n-1$ derivatives as it will be not
possible to distinguish who's who for a large value of $n$. But, again, 
read the code and find a remark with the keyword ``further derivatives'' to
see how to implement more derivatives.

\end{enumerate}

\section[Final remarks] 
		{Final remarks} 
\label{Conclusions}

\par \hskip 0.5cm As  $f=\chi^{n}$ is a {\em kernel}, in the sense that it is
positive and its integral is one, then it is integral preserving under
convolution product. Well, this is exactly the reason why $f$ has
this property because $\chi$ is a {\em kernel}. This permits to construct
a partition of the unity from a set of characteristic functions
associated with a partition of an interval with the desired class of
differentiability. This can be done by making the convolution of $f$ with
each member of this set of characteristic functions and the methodology
of sections 2,3,4 of paper can be of help to produce this rapidly.

Given a {\em kernel}, in the above sense, we can transform it in a sequence
which is an {\em approximate identity}, this is the terminology of
\cite[chapter 6, definition 6.31]{RudinFA}, the equation is
$
\phi_{n}=nf(nx)
$.
By the way, this is the expression to have a sharp {\em kernel} when constructing
a partition of the unity, instead of using $f$ directly. A simple translation
will produce a {\em kernel} balanced around zero, and sometimes this is interesting to have. The program presented in this paper can be easily modified to make these
applications. 

A simple consequence of theorem (\ref{TeoremaDois}), applied to the an $n$-splines,  show that 
its derivative of order  $n-1$ will be a linear combination of characteristic functions
of the elements of a partition of an interval of $\mathbf R$ (which can be the whole line).
This permits a simple characterization of a quite general class of univariated splines.

\begin{theorem}[Characterization]{spaces of splines}
\label{TeoremaCinco}

Given a partition $\Pi(I)$ of an interval $I$ of $\mathbf R$ (which can be the whole line), and a {\em kernel} $\phi$ being an $n$-splines of compact support
\begin{enumerate}
\item let's consider $\Pi(I)$ as the index set of the characteristic functions
of the intervals in $\Pi(I)$;
\item the convolutions of $\phi$ with the elements of $\Pi(I)$ are linearly 
independent and generate a vector space of functions of class ${\cal C}^{n-1}$
whose dimension is card($\Pi(I)$),  $Spl_{\Pi,\phi}(I)$;
\item the elements of $Spl_{\Pi,\phi}(I)$ is a space of $n$-splines.
\end{enumerate}

\end{theorem}

Observe that that there is a great generality of objects which can be obtained
this way, this diversity resulting from the possible shapes of $\phi$. 
In case $\phi$ is
the power of a characteristic function, the methodology of this paper applies to the
construction of $Spl_{\Pi,\phi}(I)$. You can have a glimpse of the power of this
methodology looking at the graphic of [a bad example]\cite{Praciano:Um}.

\bibliographystyle{unsrt}  
\bibliography{tarcisioblb}{}  
\end{document}